\documentclass[11pt]{amsart}
\usepackage{amscd}
\usepackage{amssymb}
\input xy
\xyoption{all}
\newdir{ >}{{}*!/-5pt/@{>}}
\newcommand{\ay}{\mathfrak{A}}
\newcommand{\ca}{\mathfrak{C}}
\newcommand{\ij}{\mathfrak{I}}
\newcommand{\B}{\mathfrak{B}}
\newcommand{\barijo}{\tilde{\ij}}

\newcommand{\cha}{{\rm{Ch}}{\ay}}

\newcommand{\wiso}{\overset{\cong}{\longrightarrow}}

\newcommand{\cofi}{\rightarrowtail}
\newcommand{\fib}{\twoheadrightarrow}
\newcommand{\defor}{\overset\sim\fib}
\newcommand{\fake}{``\prod\text{''}}
\DeclareMathOperator{\ch}{Ch}
\DeclareMathOperator{\Hom}{HOM}

\DeclareMathOperator{\coker}{coker}
\DeclareMathOperator{\pro}{pro}
\DeclareMathOperator*{\colim}{colim}
\DeclareMathOperator{\alg}{Alg}

\newcommand{\C}{\mathbb{C}}
\newcommand{\Z}{\mathbb{Z}}
\newcommand{\N}{\mathbb{N}}
\newcommand{\Q}{\mathbb{Q}}
\newcommand{\HP}{\mathbb{HP}}

\newcommand{\F}{\mathcal{F}}
\newcommand{\G}{\mathcal{G}}

\newcommand{\Ca}{\mathcal{C}}
\newcommand{\Tau}{\mathcal{T}}

\theoremstyle{plain}
\newtheorem{thm}{Theorem}[subsection]

\newtheorem{lem}[thm]{Lemma}
\newtheorem{prop}[thm]{Proposition}

\newtheorem{them}{Theorem}[section]
\newtheorem{coro}[them]{Corollary}
\newtheorem{lemm}[them]{Lemma}
\newtheorem{propo}[them]{Proposition}

\theoremstyle{definition}

\newtheorem{rem}[thm]{Remark}

\newtheorem{defi}[them]{Definition}
\newtheorem{rema}[them]{Remark}
\newtheorem{nota}[them]{Notation}
\newtheorem{exam}[them]{Example}

\theoremstyle{remark}

\newtheorem*{ack}{Acknowledgement}
\newtheorem*{stass}{STANDING ASSUMPTION}

\begin{document}

\title[Excision in bivariant periodic cyclic cohomology]{Excision in bivariant periodic 
cyclic cohomology: a categorical approach}

\author[G. Corti\~nas]{Guillermo Corti\~nas*}

\address{Departamento de Matem\'atica\\
Ciudad Universitaria\\
 Pabell\'on 1\\
(1428) Buenos Aires\\
 Argentina.}
\email{gcorti@dm.uba.ar}
\author[C. Valqui]{Christian Valqui**}
               
\address{IMCA\\
         Jr. Ancash 536, Casa de las Trece Monedas\\
         Lima 1, Per\'u.}
\email{cvalqui@pucp.edu.pe}
\begin{abstract}
We extend Cuntz-Quillen's excision theorem for algebras and pro-algebras in arbitrary $\Q$-linear
categories with tensor product. 
\end{abstract}

\thanks{(*) CONICET Researcher and ICTP Associate. Partially supported by grant
UBACyT X066.\\
(**) IMCA-PUCP; partially supported by CONCYTEC grant CS017-2002-OAJ}

\maketitle

\section{Introduction}
Bivariant periodic cyclic cohomology is a bifunctor defined on the category of algebras
over a field, contravariant in the first variable and covariant in the second, which associates
a $\Z/2$-graded vectorspace $HP(A,B)$ with every pair of algebras $(A,B)$. Cuntz-Quillen excision
theorem (\cite{cq}) says that, in characteristic zero, to every exact sequence of algebras
\begin{equation}\label{sexseq}
0\to A'\to A\to A''\to 0
\end{equation}
and every algebra $B$ there correspond natural homogeneous maps 
$$
HP(A',B)\to HP(A'',B)[1]\text{ and } HP(B,A'')\to HP(B,A')[1]
$$
such that the triangles of $\Z/2$-graded vectorspaces
\begin{align}\label{triangles}
\xymatrix{HP(A'',B)\ar[rr]&& HP(A,B)\ar[dl]\\
            &HP(A',B)\ar[ul]&}\\
\xymatrix{HP(B,A')\ar[rr]&& HP(B,A)\ar[dl]\\
            &HP(B,A'')\ar[ul]&}
\end{align}
are exact. This result was extended first 
to a special type (\cite{cu}) and then to all 
 topological algebras with jointly continuous multiplication (\cite{v1}).
The latter are complete locally convex $\C$-vectorspaces $A$ equipped with an associative 
multiplication $A\otimes_\pi A\to A$, where $\otimes_\pi$ is the projective tensor product.
Extensions of Cuntz-Quillen's result in other directions were given in \cite{ralf} and \cite{pusch}. 
For example in \cite{ralf} it was shown --among other things-- that excision holds for $HP$ of discrete 
pro-algebras. These are inverse systems $\{A_n\to A_{n-1}\}_{n\in\N}$ of (discrete) algebras. 
In both the topological and the pro-setting excision holds under the assumption that
\eqref{sexseq} be split in an underlying category; that of locally convex
vectorspaces in the first case, and that of algebraic pro-vectorspaces in the second. 
Here we further extend the excision
theorem to algebras and pro-algebras in $\Q$-linear categories (additive categories with uniquely
divisible $\hom$ groups) with tensor product, as follows. 
Let $\ay$ be a $\Q$-linear category closed under finite limits and colimits and 
$\otimes:\ay\times\ay\to \ay$ a symmetric, strictly associative product, such that for $V\in\ay$,
the functor $V\otimes$ preserves cokernels and split exact 
sequences. An algebra
in $\ay$ is an object together with an associative ``multiplication'' map $A\otimes A\to A$. Write
$\alg\ay$ for the category of algebras. We define bivariant periodic cyclic cohomology for pairs of algebras in $\ay$ just as for concrete algebras; we put
\begin{equation}\label{hp}
HP^*(A,B):=H^*\Hom_{\Z/2}(X^\infty A,X^\infty B)
\end{equation}
Here $X^\infty$ is the Cuntz-Quillen pro-supercomplex (\cite{cq}), and $\Hom_{\Z/2}(,)$ is the mapping 
supercomplex.
We show that if $B, A', A, A''$ are objects of
$\alg\ay$ and  
\eqref{sexseq} is split exact then the triangles \eqref{triangles} are exact. 
The same holds for $B,A',A,A''\in\pro$-$\alg\ay$, as pro-algebras form a subcategory of the category
of algebras in $\ay'=\pro$-$\ay$, and $\ay'$ is $\Q$-linear. 
However we show that it is also possible to
extend periodic cyclic cohomology to pro-algebras in a different way, and thus obtain exact triangles 
as in \eqref{triangles} under a weaker assumption than that \eqref{sexseq} be split as sequence 
in $\pro$-$\ay$. Note
\eqref{sexseq} is split in $\pro$-$\ay$ if and only if every map $A'\to V\in\pro$-$\ay$
can be extended to a map $A\to V$; pictorially

\begin{equation}\label{rlp}
\xymatrix{A'\ar[r]\ar@{ >->}[d]& V\\
          A\ar@{.>}[ur]}
\end{equation}
Here we prove that excision holds whenever this extension
property is fulfilled just for objects $V\in\ay$ (viewed as constant pro-objects) rather than for all 
pro-objects. This includes the case when $A'\to A$ can be represented 
by a map of inverse systems $i_n:A'_n\to A_n$ such that each $i_n$ is individually split but
the splittings do not commute with the transfer maps $A_n\to A_{n-1}$, $A'_n\to A'_{n-1}$ 
(cf. \ref{chri}).
To obtain this stronger excision result we replace periodic cohomology \eqref{hp} by hypercohomology; 
we prove excision for 
$$
\HP^*(A,B):=H^*\Hom_{\Z/2}(X^\infty A,\widehat{X}^\infty B)
$$
substituted for $HP$ (see \ref{6tes}). Here the $\widehat{\ \ }$ denotes any fibrant resolution of $X^\infty B$ in a certain 
closed model structure for pro-supercomplexes which we introduce in section \ref{seccm}. If $B$
is an algebra then $X^\infty B$ is fibrant already and thus
$$
\HP^*(A,B)=HP^*(A,B)
$$
in this case.

We remark that the categories of algebraic, complete locally convex and bornological $\C$-vectorspaces 
are $\Q$-linear categories with tensor product in our sense. Thus excision for 
$HP$ of algebraic, topological and bornological algebras and pro-algebras follows from our results.

An interesting feature of our proof is that we do not assume infinite sums exist in 
$\ay$. In particular we have no analogue for the tensor algebra of a vectorspace, much
less for a left adjoint of the forgetful functor $\alg\ay\to \ay$. This
is remarkable, because the lack of an obvious way to topologize the tensor algebra
of a complete locally convex vectorspace in such a way as to have good adjointness properties
--later found in \cite{v1}-- initially restricted the proof of excision to topological algebras 
whose topology is defined by a family of submultiplicative seminorms (\cite{cu}). 
Our current methods evolved from, strengthen and generalize those of \cite{v2}.

\bigskip

The rest of this paper is organized as follows. In section \ref{cplx} we consider pairs
$(\ay,\ij)$ where $\ay$ is an additive category and $\ij$ any class of objects of $\ay$,
and introduce the term {\it cofibration} for maps $f:A'\cofi A\in\ay$ 
having the extension property \eqref{rlp} for every $V\in\ij$. We call an object $V\in\ay$
{\it relatively injective} if the extension property \eqref{rlp} holds whenever $A'\cofi A$
is a cofibration. For example objects of $\ij$ are relatively injective, but there may be more.
Note that to say that every object is relatively injective is the same as to say every 
cofibration is a split injection. Further if $A'\to A$ is a map in $\pro$-$\ay$, we call it
a cofibration if it has \eqref{rlp} for all relatively injective $V\in\ay$.
Starting with these basic notions, we introduce 
Waldhausen category structures (\cite{wald}) for the various categories of complexes and pro-complexes 
relevant to  cyclic homology, and show that this added structure is preserved by the usual 
functors going between them. In section \ref{seccm} we show that if every
object of $\ay$ (but not $\pro$-$\ay$) is relatively injective, then the notions of 
cofibration and weak equivalence for pro-supercomplexes in $\ay$ introduced on the previous
section, together with an appropriate class of fibrations, make $\pro$-$Sup\ay$ into a
closed model category in the sense of Quillen (\cite{qui}).
In section \ref{alg} we work out the construction of 
the relative complexes 
relevant to Wodzicki's excision theorem for our abstract algebras $A\in\alg\ay$.
Section \ref{mariusz} is devoted to the proof of Wodzicki's excision theorem for abstract
algebras. The theorem is stated in terms of weak equivalences of complexes, as defined
in section \ref{cplx}. The proof works for arbitrary pairs $(\ay,\ij)$ as above; in particular
cofibrations in $\ay$ need not split. In the remaining three sections we concentrate
on the case when cofibrations are split in $\ay$ (but not in
$\pro$-$\ay$).
In section \ref{exkinfi} we prove excision for pro-ideals of the form $K^\infty$ where $K$ is
 an ideal of a quasi-free (abstract) pro-algebra. Also
in this section the analogue of the Cuntz-Quillen pro-algebra $TA/JA^\infty$ is introduced;
this is always defined, even if $TA$ is not. The next section is devoted to proving Goodwillie's
theorem for abstract pro-algebras; this is the place where $\Q$- rather than $\Z$-linearity
is needed. Finally in section \ref{seccq} we prove the main theorem (Thm. 8.1; see also
Cor. 8.4).
 
\section{Complexes}\label{cplx}
\subsection{Cofibrations}\label{duno}
Throughout this section, $\ay$ will be a fixed additive category, closed under finite
limits and finite colimits.
We also fix a class $\ij$ of objects of $\ay$, which we call {\it basic}. The 
pair $(\ay,\ij)$ is a {\it basic pair}. A map 
$i:A'\cofi  A$ is a {\it cofibration} if $\hom(i,I)$ is surjective 
for all $I\in\ij$.
One checks that
$\ay$ together with the cofibrations just defined,
is a category with cofibrations in the sense of Waldhausen (\cite{wald}). 
We borrow from {\it loc. cit.} the
standard notation for cofibration sequences; in particular, if $A\cofi B$ is a cofibration
we write $B\to B/A$ for its cokernel and call
\begin{equation}\label{coseq}
A\cofi B\to B/A
\end{equation}
a {\it cofibration sequence}. 
An object $F\in\ay$ is {\it relatively injective} if $\hom(,F)$ maps cofibrations to
surjections. Equivalently $F$ is relatively injective if for every cofibration sequence
\eqref{coseq} the sequence of abelian groups
\begin{equation}\label{seq}
0\to \hom(B/A, F)\to \hom(B,F)\to\hom(A,F)\to 0
\end{equation}
is exact. For example all basic objects are relatively injective. We write $\barijo$ for
the class of relatively injective objects.
Consider the category $\cha$ of chain complexes. If $A,B\in\cha$ we put $\hom(A,B)$ for the 
homomorphisms in $\cha$ and $\Hom(A,B)$ for the chain complex whose term of degree $n$ is the 
set of homogeneous maps of degree $n$, and whose boundary operator is 
$$
f\mapsto [\partial,f]:=\partial^Bf-(-1)^{|f|}f\partial^A
$$
We identify $\ay$ with the full subcategory
of all those $A\in\cha$ with $A_n=0$ for all $n\ne 0$. A map $i:A\to B$ in $\cha$ is 
a cofibration if $\Hom(i,I)$ is degreewise surjective for all $I\in\ij$. Equivalently, $i$ is a 
cofibration if $i_n:A_n{\to}  B_n$ is one for all $n$. For example 
if $f:A\to B$ is a chain map and $C_f$ is its mapping cone, then the inclusion 
$B\subset C_f$
is a cofibration (because it is a split mono degreewise) and 
\begin{equation}\label{mc}
0\to B\to C_f\to A[-1]\to 0
\end{equation}
is a cofibration sequence.
An object $F\in\cha$ is called {\it fibrant}
if $F_n$ is relatively injective for all $n$. 
A chain map $f:A\to B$ is a {\it weak equivalence} if $\Hom(f,I)$ is a quasi-isomorphism
or quism for every relatively injective $I\in\ay$. One checks that with the 
classes of cofibrations and of weak equivalences 
just defined, $\cha$ is a category with cofibrations and weak equivalences in the sense of \cite{wald}
satisfying both the saturation and the extension axioms of {\it loc. cit.} 
An object $Q\in\cha$ is {\it weakly contractible} if $0\to Q$
--or equivalently $Q\to 0$-- is a weak equivalence.

\begin{lem}\label{weq}
Let $f:A\to B\in\cha$, $C$ its mapping cone. The following conditions are equivalent
for $f:A\to B\in\cha$. 
\item{i)}  $f$ is a weak equivalence.
\item{ii)}  $C$ is weakly contractible.
\item{iii)}  For every $n$, the map
$$
\frac{C_{n}}{\partial C_{n+1}}\overset{j}{\cofi} C_{n-1}
$$
induced by $\partial$ is a cofibration.

\item{iv)}  For every fibrant complex $F$ which is bounded (both above and below), 
the map $\Hom(f,F)$ is a quism.

\end{lem} 
\begin{proof}
That (i)$\iff$(ii) is immediate from the fact that \eqref{mc} is a cofibration sequence. 
That (ii)$\iff$(iii) is
immediate from the definition of weakly contractible and the left exactness of $\hom$. It is clear 
 that (iv)$\Rightarrow$(i). To prove the converse, note first that, because the sequence \eqref{mc}
is degreewise split, (iv) is equivalent to saying that $\Hom(C,F)$ is exact for all bounded fibrant $F$. 
We know this holds for $F$ of length one, since (i)$\Rightarrow$(ii). For $F$ of 
length $n\ge 2$ it follows by induction as we apply $\Hom(C,)$ to the degreewise split exact sequence
$$
0\to F_{\le n-1}\to F\to F_{ n}\to 0
$$
\end{proof}


\subsection{Pro-objects}\label{pro}
If $\ca$ is a category, we write $\pro$-$\ca$ for the category whose objects are the inverse systems
$$
A=\{A_{n+1}\overset{\sigma}{\to} A_{n}: n\ge 1\}
$$
indexed by $\N:=\Z_{\ge 1}$, and where the homomorphisms are the pro-maps
$$
\hom(A,B):=\lim_{n}\colim_m\hom(A_m,B_n)
$$
We note our definition differs from that of \cite{cq}, but yields an equivalent category. 
It is useful to think of pro-homomorphisms as equivalence classes of maps of inverse systems.
For this purpose, we introduce some notation. If $A,B\in\pro$-$\ca$, then a {\it representative map} 
from $A$ to $B$ is a nondecreasing, nonstationary function $f:\N\to \N$
together with a map of inverse systems (which we also denote by $f$)
$$
f:f_*A\to B
$$
Here
$$
f_*A_n:=A_{f(n)}
$$
If $f,g$ are representatives we put $f\le g$ if first of all 
$f(n)\le g(n)$, and second of all
\begin{equation}
\xymatrix{A_{g(n)}\ar[r]^g\ar[d]_\sigma &B_n\\
          A_{f(n)}\ar[ur]_f&}
\end{equation}

commutes, for every $n\in\N$. The relation $\le$ is a partial order in the set of all representatives 
from $A$ to $B$. We say that two representatives $f$ and $g$ are {\it equivalent} if they have a common
upper bound. 

\begin{lem}\label{rep}
Let $\ca$ be a category, $A,B\in\pro$-$\ca$, $R(A,B)$ the set of all representatives from $A$ to 
$B$ and $\equiv$ the equivalence relation defined above. Consider the map
$$
\gamma:R(A,B)\to \hom(A,B),\qquad f\mapsto ([f_n])_{n\in\N}
$$
where $f_n:A_{f(n)}\to B_n$ is the $n$-th component of $f$ and $[f_n]$ its class in 
$\colim_m\hom(A_m,B_n)$. Then $\gamma$ induces a bijection
$$
R(A,B)/\equiv\overset{\cong}\longrightarrow \hom(A,B)
$$
\end{lem}
\begin{proof} Straightforward.
\end{proof}

Note that $\ca$ is isomorphic to the full subcategory of $\pro$-$\ca$ consisting of those
pro-objects such that all structure maps $\sigma$ are identity maps ({\it constant} pro-objects).
We shall identify $\ca$ with this subcategory of $\pro$-$\ca$. If $\ay$ is an additive category
in which a class $\ij$ of basic objects has been chosen, then we choose $\barijo$ as the
basic objects of $\pro$-$\ay$. To give examples and --under extra hypothesis-- obtain a full 
characterization of the relatively injective pro-objects which result from this choice, 
we need to recall the
following well-known construction.
The {\it fake product} of a sequence $A=(A_n)_{n}$ of 
objects of $\ay$ is the pro-object $``\prod$''$A$ given by
\begin{equation}\label{pra}
{\fake}_n A
=\oplus_{p=0}^nA_p
\end{equation}
with the obvious projections as structure maps. If $X\in\pro$-$\ay$, then
\begin{equation}\label{hompra}
\hom(X, ``{\prod}\text{''}A)=\prod_{n=1}^\infty\hom(X,A_n)
\end{equation}
Now suppose $A$ is not just a sequence of objects but a pro-object. Then for each $n$ the
identity map $1_n$ of $A_n$ represents a map $A\to A_n$. Putting all these together and using
\eqref{hompra} we get a map 
\begin{equation}\label{iota}
\iota:A\cofi \fake A
\end{equation}
Note that $\iota$ is a cofibration. Indeed, for $I\in\ay$, 
$$\hom(\fake A,I)
=\bigoplus_{n=1}^\infty\hom(A_n,I)\twoheadrightarrow \colim_n\hom(A_n,I)=\hom(A,I)$$
is the natural projection. 
\begin{lem}\label{equi-inj}
Let $A\in\pro$-$\ay$. Consider the following conditions.
\item{i)} There exists a sequence of relatively injective objects of $\ay$, $X=(X_n)_n$, such that
$A$ is a retract of $``\prod$''$X$.

\item{ii)} $A$ is relatively injective. 

Then i)$\Rightarrow$ii). The converse holds under the following extra assumption:

\noindent{For every $B\in\ay$ there exists a cofibration $B\cofi I$ with $I$
 relatively injective.}

\end{lem}
\begin{proof}
A retract of a relatively injective object is relatively injective. Thus to prove that i)$\Rightarrow$ii)
it suffices to show that the fake product of a sequence of relatively injective objects is relatively
injective. But this is immediate from \eqref{hompra}, since $\hom(,X_n)$ maps cofibrations to surjections.
Next assume the extra hypothesis holds and
let $A$ be relatively injective. Choose a cofibration $j_n:A_n\cofi I_n$ for each $n\ge 1$. Each $j_n$
can be viewed as element of $\oplus_m\hom(A_m,I_n)=\hom($``$\prod$''$A,I_n)$, and the collection 
$j$ of all the $j_n$ as an element of 
$\prod_n\hom($``$\prod$''$A,I_n)=\hom($``$\prod A$'',``$\prod$''$I)$ (using \ref{hompra}).
If $J$ is a constant pro-object, then $\hom(j,J)$ is the map
$$
\xymatrix{
\bigoplus_n\hom(j_n,J):&\bigoplus_n\hom(I_n,J)\ar@{>>}[r]\ar@{=}[d]&\bigoplus_n\hom(A_n,J)\ar@{=}[d]\\
&\hom(\fake I,J) &\hom(\fake A,J)}
$$
Thus $j$ is a cofibration. Now consider the cofibration $\iota$ of \eqref{iota}. Because $A$ is relatively
injective, $j\iota$ is a split injection. Hence $A$ is a retract of ``$\prod$''$I$.
\end{proof}

The extra condition of the lemma is fulfilled for example if every object of $\ay$ is relatively
injective. Thus in this case $A\in\pro$-$\ay$ is relatively injective $\iff$ $\iota:A\to``\prod$''$A$ is
a split injection.  Next we characterize for $(\ay,\ij)$ arbitrary, the pro-objects $A$ such that 
$\iota$ is a split injection.
Let $A=\{\sigma_n:$ $A_{n}\to A_{n-1}\}\in\pro$-$\ay$; we say that $A$ has property $(P)$ if there
exists a strictly increasing sequence $\{n_k:k\ge 1\}$ such that
$$
(P)\quad \forall k\ \ \exists s_k:A_{n_k}\to A_{n_{k+1}}\text{ such that } \sigma^{n_{k+1}-n_{k-1}}
s_k=\sigma^{n_k-n_{k-1}}
$$
Here $\sigma^p:A_{*+p}\to A_*$ is the structure map. Note that  $A$ is isomorphic to $A':=\{A_{n_k}\}_k$.
Thus we can always replace $A$ by an isomorphic pro-object satisfying $(P)$ for $n_k=k+1$ ($k\ge 1$).

\begin{lem}
Let $A\in\pro$-$\ay$, $\iota:A\cofi``\prod$''$A$ the canonical cofibration. The following are equivalent
\item{i)} $\iota$ is a split injection.
\item{ii)} $A$ has property $(P)$.
\end{lem}
\begin{proof} i)$\Rightarrow$ii): Since $``\prod$''$A$ always has $(P)$, it suffices to show
that if $B$ has $(P)$ and $i:A\to B$ is a split injection then $A$ does too. Let $p:B\to A$
be a left inverse for $i$. Without loss of generality we may assume $p,i$ have representatives
$\{i:A_n\to B_n\}_n$, $\{p:B_{n+1}\to A_n\}_n$ such that $pi=\sigma$, and that there is a
sequence of maps $s_n:B_n\to B_{n+1}$ with $\sigma^2 s_n=\sigma$. Put $t=ps^2i:A_n\to A_{n+1}$.
One checks that 
$$\sigma^3t=\sigma^2 $$
whence $\sigma^4t^2=\sigma^2$. Thus $A$ verifies property $(P)$ with $n_k=2k-1$ ($k\ge 1$).

ii)$\Rightarrow$i) Let $(n_k)_{k\ge 1}$ and $(s_k)_{k\ge 1}$ be as in $(P)$. Consider $f:\N\to \N$, $f(k)=n_k$.
We have a commutative
diagram 
$$
\xymatrix{f_*A\ar[r]^\cong\ar[dr]_{\iota_{f_*A}}&A\ar[d]^p\ar[r]^{\iota_A}&\fake A\ar[dl]^q\\
           &\fake f_*A&}
$$
Here $p\in\hom(A,``\prod$''$f_*A)=\prod\hom(A,A_{n_k})$ is given by the canonical maps $A\to A_{n_k}$
$(k\ge 1)$ and $q$ by composite maps 
$$
\fake A\to{\fake}_{n_k}A\fib A_{n_k}
$$
It follows that if $\iota_{f_*A}$ is a split monomorphism, then so is $\iota_A$. Therefore
we may assume $n_k=k$ for all $k\ge 1$. Put
$$
r_n:{\fake}_{n+2}A=\bigoplus_{p=1}^{n+2}A_p\to A_n
$$
for the map given by the transpose of the following matrix of homomorphisms
$$
[0,\sigma s^n\sigma, \sigma s^{n-1}(\sigma-s\sigma^2),\dots,\sigma s(\sigma-s\sigma^2),\sigma
(\sigma-s\sigma^2)]
$$
Here we have omitted the subscripts of $s$. One checks that $r_n\sigma=\sigma r_{n+1}$ and 
$r_n\iota_{n+2}=\sigma^2$ ($n\ge 1$).
\end{proof}

In the case when all objects of $\ay$ are relatively injective, that is, when all cofibrations
in $\ay$ are split monomorphisms, we have the following simple
characterisation of cofibrations in $\pro-\ay$.

\begin{lem}\label{chri}
Let $f:X\to Y$ be a map in $\pro-\ay$. Assume all objects of $\ay$ are relatively injective.
Then $f$ is a cofibration if and only if there exists a commutative diagram
$$
\xymatrix{X\ar[r]^f\ar[d]^\wr &Y\ar[d]^\wr\\
 X'\ar[r]_{f'}&Y'}
$$
where the vertical maps are isomorphisms and  $f'$ has a representative map $\{f'_n:X'_n\cofi Y'_n\}_n$ such that each $f'_n$ is
a split monomorphism. 
\end{lem}

\begin{proof}
Without loss of generality we may assume $f$ is represented by a level map $\{f_n:X_n\to Y_n\}_n$.
By \ref{equi-inj}, $f$ is a cofibration $\iff$ the dotted arrow in the diagram below
exists
$$
\xymatrix{X\ar[r]^(.4){i}\ar@{ >-}[d]_f&\fake X\\
          Y\ar@{.>}[ur]_s&}
$$
One checks (using \eqref{hompra}) that this condition is equivalent to the existence, for each
$k\ge 1$, of an $n_k\ge 1$ and a map $s_k:Y_{n_k}\to X_k$ such that 
\begin{equation}\label{deco}
s_kf_{n_k}=\sigma^{n_k-k}
\end{equation}
This condition is certainly satisfied by the map $f'$ of the lemma; since in addition 
cofibrations as defined above are clearly preserved by isomorphisms of arrows, it follows
that the second condition of the lemma implies the first. Next suppose $f$ is a cofibration;
upon replacing it by an isomorphic arrow, we may assume $n_k=k+1$ in \eqref{deco}. Put
$$
\sigma':=\left[\begin{matrix} 0&s\\
                       0 &\sigma\end{matrix}\right]\qquad j=(s,\sigma)
$$
We have a commutative diagram
\begin{equation}\label{mispi}
\xymatrix{X_{n+1}\ar[d]_\sigma\ar[r]^(.4){(1,f)}& X_{n+1}\oplus Y_{n+1}\ar[d]^{\sigma'}\ar[r]^(.6){\pi}&
Y_{n+1}\ar[dl]^{j}\ar[d]^\sigma\\
X_n\ar[r]_(.4){(1,f)} & X_n\oplus Y_n\ar[r]_(.6){\pi}&Y_n}
\end{equation}
Here $\pi$ is the obvious projection.
Write $Y'=\{X_n\oplus Y_n,\sigma'\}_n$. Using \eqref{mispi}, we obtain a pro-map
$f':X\to Y'$ and a pro-isomorphism $j:Y\to Y'$ (with inverse $\pi$) which fit into
a commutative diagram 
$$
\xymatrix{X\ar@{=}[d]\ar[r]^f&Y\ar[d]^{\wr}_j\\
          X\ar[r]_{f'}&Y'}
$$
\end{proof}

\subsection{}\label{procha}
Because $(\pro$-$\ay,\barijo)$ is a basic pair, in  $\ch\pro$-$\ay$ there are defined notions of 
cofibration, fibrant object and weak equivalence. On the other hand we may also 
consider the category $\pro$-$\cha$. There is an obvious functor
$\gamma:\pro$-$\cha\to\ch\pro$-$\ay$ induced by the natural inclusion of the class of objects of
the first category into the second, which is neither full nor faithful. We say that a map 
$f$ in $\pro$-$\cha$ is a cofibration (resp. a weak equivalence) if $\gamma(f)$ is, and that
an object $F$ is fibrant if so is $\gamma(F)$. One checks that $\gamma$ preserves pushouts,
whence it follows that the cofibrations in $\pro$-$\cha$ satisfy Waldhausen's axioms, as do
the weak equivalences. In other words, $\gamma$ is an {\it exact functor} between
Waldhausen categories (\cite{wald}).

\subsection{Mixed complexes}\label{mix}
We adopt the definitions of \cite{cq2}, except that we work
with complexes of objects of $\ay$ rather than $\C$-vectorspaces. In particular, mixed 
complexes are assumed to be bounded below. We write $M\ay$
for the category of mixed complexes.

A map in $M\ay$ is a {\it cofibration} (resp. a {\it weak equivalence}, resp. an object is {\it fibrant}) 
if it is mapped
to one by the forgetful functor $(M,b,B)\mapsto (M,b)$. One checks that this structure makes $M\ay$ into
a Waldhausen category for which the forgetful functor is exact.
If $(M,b,B)$ and $(N,b,B)$ are mixed 
complexes then the chain complex $\Hom((M,b),(N,b))$ is a mixed complex with 
Connes operator $[B,]$. It is clear that a map $f\in M\ay$ is a cofibration (resp. a weak equivalence)
if for the mapping mixed complex, $\Hom(f,I)$ is surjective (resp. is a quism of the underlying chain
complexes) for every $I\in\barijo$. Similarly an object $F$ is fibrant $\iff$ $\Hom(,F)$ sends cofibrations
to surjections. One checks that if $f:M\to N$ is a map of mixed complexes,
then the mapping cone $P$ of $f$ as a map in $\cha$ is again a mixed complex. 
The Connes
operator is defined as
$$
\left(\begin{matrix} B &0\\
         0&-B\end{matrix}\right):P_n=N_n\oplus M_{n-1}\to P_{n+1}
$$
Parts i)$\iff$ ii) and ii)$\iff$iv) of Lemma \ref{weq} hold for mixed complexes,
with the mapping cone just defined. 

\subsection{$S$-complexes}\label{Scom}
We write $S$-$\ay$ for the category of bounded below chain complexes of objects of $\ay$
equipped with a periodicity operator $S:P\to P[-2]$. 
We give $S$-$\ay$ the only Waldhausen category structure which makes the forgetful functor
$(P,\partial,S)\mapsto (P,\partial)$ exact. 
The mapping
$S$-complex of two objects $P,Q\in S$-$\ay$ is defined in terms of the ordinary mapping chain
complex as
$$
\CD
\Hom_S(P,Q):=\ker(\Hom(P,Q)@>[S,]>>\Hom(P,Q)) 
\endCD
$$
It is possible to express the Waldhausen structure of $S$-$\ay$ in terms of this mapping complex,
as follows. A map $f\in S$-$\ay$ is a {\it cofibration} (resp. a {\it weak equivalence}) $\iff$ for each $I\in\barijo$
--considered as an object of $M\ay$-- $\Hom(f,\B I)$ is surjective (resp. a quism).
Here $\B:M\ay\to S$-$\ay$ is the usual {\it bar construction}
\begin{equation}
\B M_n=\bigoplus_{p\ge 0}M_{n-2p}
\end{equation}
The equivalence between the original notions of cofibration and weak equivalence and the
ones given in terms of $\Hom_S$ is immediate from the particular case of the identity
\begin{equation}\label{hom(,B)}
\Hom_S(P,\B N)_n=\prod_r\hom(P_r,N_{n+r})
\end{equation}
when $N$ is concentrated in degree zero.
We say that $F\in S$-$\ay$ is {\it fibrant} if $\Hom_S(,F)$ maps cofibrations to surjections.
For example --as follows from \eqref{hom(,B)}-- $\B M$ is fibrant if $M$ is a fibrant mixed
complex.
 
If $f:P\to Q$ is a map in $S$-$\ay$, then the ordinary mapping cone of $f$ as a chain map
is again an $S$-complex, with the periodicity operator induced by those of $P$ and $Q$.
Part i)$\iff$ii) of Lemma \ref{weq} holds in $S$-$\ay$ for this mapping cone; the argument of the proof
of part i)$\iff$iv) shows that a map $f:P\to Q\in S$-$\ay$ is a weak equivalence $\iff$ 
$\Hom(f,\B F)$ is a quism for every bounded fibrant mixed complex $F$.

\begin{prop}\label{Bexact}
The bar construction $\B:M\ay\to S$-$\ay$ is an exact functor of Waldhausen categories which preserves
fibrant objects.
\end{prop}
\begin{proof}
In both $M\ay$ and $S$-$\ay$, a map is a cofibration $\iff$ it is so degreewise. Thus if 
$f$ is a cofibration then $\B f$ is, because in each degree it is a finite direct sum of cofibrations.
That $\B$ preserves fibrant objects is immediate from \eqref{hom(,B)}; that it commutes with finite
pushouts is clear. To show that also weak equivalences are preserved, it suffices to show $\B M$ is
weakly contractible if $M$ is. But for $I\in\ay$
\begin{equation}\label{hom(B,I)}
\Hom_S(\B M,\B I)_n=\bigoplus_{p\ge 0}\hom(M_{-n-2p},I)\qquad (n\in\Z)
\end{equation}
is the direct sum total complex of the following bicomplex of abelian groups
$$
D_{p,q}=\hom (M_{q-p},I)
$$
Here the horizontal and vertical boundary maps are respectively 
$\hom(B ,I)$ and $\hom(b,I)$. 
If $M$ is weakly contractible and $I\in\barijo$, then each column of $D$ is exact, whence so is
$\Hom_S(\B M,\B I)$.
\end{proof} 

\subsection{Supercomplexes, pro-supercomplexes and
 super-pro-\\
complexes}\label{sup}
Recall a supercomplex in $\ay$ is a $\Z/2$-graded object with a square-zero differential of degree $1$.
We write $Sup\ay$ for the category of supercomplexes. 
If $X,Y\in Sup\ay$ then $\Hom_{\Z/2}(X,Y)$ is the supercomplex of $\Z/2$-homogeneous maps.
We identify $Sup\ay$ with the full subcategory of $Sup\pro$-$\ay$ of supercomplexes of constant pro-objects.
We say that a map $f\in Sup\pro$-$\ay$ is a {\it cofibration } (resp. a {\it weak equivalence }) if
$\Hom_{\Z/2}(f,Y)$ is surjective (resp. a quism) whenever $Y\in Sup\ay$ and $Y_i$ is a
relatively injective object of $\ay$ ($i=0,1$). This structure makes $Sup\pro$-$\ay$ into a Waldhausen
category.
An object $X\in Sup\pro$-$\ay$ is {\it fibrant } if $\Hom_{\Z/2}(,X)$ maps cofibrations
to surjections, or equivalently if both its even and odd degrees are relatively injective objects
of $\pro$-$\ay$. All this structure can be translated to $\pro$-$Sup\ay$ using the canonical functor
$\pro$-$Sup\ay\to Sup\pro$-$\ay$, in the same way as was done for chain complexes (\ref{procha}).  
If $M$ is a mixed complex then we can associate with it the pro-supercomplex 
$$X^\infty M=\{X^nM\}$$ of \cite{cq}, which at level $n$ is the supercomplex 
$$
X^n=\frac{M_n}{bM_{n+1}}\oplus \bigoplus_{i=0}^{n-1}M_i
$$
with its natural even-odd grading, and boundary map $X^nM\to X^nM$ induced by $b+B$. 
We also consider the
super-pro-complex $\xi M=\{\xi^nM\}$ of \cite{cu}, which is
$$
\xi^nM=\oplus_{i=0}^nM_i
$$
at level $n$, with the obvious grading and boundary map $X^n\to X^{n-1}$ also induced by $b+B$.
One checks that the canonical functor $\pro-Sup\ay\to Sup\pro-\ay$ maps $X^\infty M$ isomorphically 
onto $\xi M$.

\begin{prop}\label{xiexact}
The functor $\xi:M\ay\to Sup\pro$-$\ay$ is an exact functor of Waldhausen categories
which preserves fibrant objects.
\end{prop}
\begin{proof} Let $M$ be a mixed complex and  $Y$ a supercomplex. We have
\begin{equation}\label{pliji}
\Hom_{\Z/2}(\xi M,Y)_m
=\bigoplus_n\hom (M_n,Y_{\overline{m+n}})\qquad (m\in\Z/2)
\end{equation}
It is clear from this identity that $\xi$ preserves cofibrations. To see that it also preserves 
weak equivalences, it suffices to show that if $M$ is weakly contractible then so is $\xi M$. 
We have to prove that, for $Y$ fibrant and constant, \eqref{pliji} is exact as a $\Z/2$ or equivalently
as a $\Z$-graded complex. Consider the subcomplex
$$
D^p:=\bigoplus_{n\le p}\hom(M_n,Y_{\overline{p+n}})\subset \Hom_{\Z/2}(\xi M,Y)_{\overline{p}}\qquad (p\ge 0)
$$
One checks that every cocycle of $\Hom_{\Z/2}(\xi M,Y)$ is also a cocycle of $D$ of some dimension.
Thus we are reduced to showing that $D$ is exact. Consider the subcomplex
$$
{D'}^p:=\bigoplus_{n\le p-1}\hom(M_n,Y_{\overline{p+n}})\subset D^p
$$
There are exact sequences
\begin{align*}
0\to D'\to D&\to \Hom (M,Y_0)\to 0\\
0\to D[-2]\to D'&\to \Hom (M,Y_1)\to 0
\end{align*}
Because $M$ is weakly contractible, $\Hom (M,Y_i)$ is exact for $i=0,1$. It follows from 
the exact sequences above that both $D$ and $D'$ are exact.
Finally if $M$ is fibrant then
each $M_n$ is relatively injective, whence $\xi_rM=$``$\prod$''$M_{2*+r}$ is relatively injective
in $\pro$-$\ay$ (by \ref{equi-inj}) and therefore $\xi M$ is fibrant.
\end{proof}

\begin{rem}\label{remxiexact}
Let $f:M\to N$ be a weak equivalence of pro-mixed complexes. This means that, for the canonical functor
$\gamma:\pro$-$M\ay\to M\pro$-$\ay$, $\gamma(f)$ is a weak equivalence. Hence by
\ref{xiexact} , $X^\infty\gamma(f)$ is a weak equivalence in $Sup\pro$-$\pro$-$\ay$. On the other
hand $f$ also induces a map 
$$
\{X^nM^n\}_n\overset\sim\to\{X^nN^n\}_n\in\pro-Sup\ay
$$
where the superscript indicates level. I claim that also this map is a weak equivalence. The claim 
follows from the fact that, for $Y\in Sup\ay$, 
\begin{align*}
\Hom_{\Z/2}(X^\infty M,Y)=&\colim_r\Hom_{\Z/2}(X^rM,Y)\\
                            =&\colim_r\colim_s\Hom_{\Z/2}(X^rM^s,Y)\\
                            =&\colim_n\Hom_{\Z/2}(X^nM^n,Y)\\
                            =&\Hom_{\Z/2}(\{X^nM^n\}_n,Y)
\end{align*}
In what follows we shall abuse notation and write $X^\infty(M)$ to mean
the pro-supercomplex $\{X^n(M^n)\}_n$.
\end{rem}

\section{A closed model structure for pro-supercomplexes}\label{seccm}
\begin{defi}
Assume all objects of $\ay$ are relatively injective. Call a map $f:X\fib Y\in\pro$-$Sup\ay$
a {\it fibration} if it is a degreewise split surjection and its kernel is fibrant. 
\end{defi}

\begin{them}\label{clmod}
Assume every object of $\ay$ is relatively injective. Then the category  $\pro$-$Sup\ay$ with cofibrations
and weak equivalences as defined in section 2 and fibrations as in the definition above, satisfies 
Quillen's axioms {\rm (\cite{qui})} for a closed
model category.
\end{them}
\begin{proof}
It is clear that the classes of cofibrations, fibrations and weak equivalences are closed 
under retracts. Thus it suffices to prove the axioms for a model category (\cite{qui}, I.0.5).
Of these, only $M1$ and $M3$ are not immediate. To prove $M1$ consider the diagram 
\eqref{M1} of \ref{M1prop} below. The case when $i$ is a weak equivalence is \ref{M1prop}.
The case when $\pi$ is a weak equivalence follows from corollary \ref{fibequiv} by a standard
argument which does not involve pro-objects (see \cite{qui}).
To prove $M3$, let $f:A\to B$ be any map in $\pro$-$Sup\ay$. Consider the composite
$$
\theta:A\to \fake A\to C\fake A
$$
of the canonical inclusions, where $C$ stands for the mapping cone of the identity. Then
$i:=(\theta,f):A\to M:=C``\prod$''$A\oplus B$ is a cofibration, the projection $p:M\to B$
is a fibration and weak equivalence, and $f=pi$. It remains to show that $f$ can also
be written as $f=qj$ with $q$ a fibration and $j$ a cofibration and weak equivalence. 
For this purpose we introduce some notation. For $X\in\pro$-$Sup\ay$ we write 
$RX$ for the mapping cone of the map of inverse systems
$$
\left[\begin{matrix} 1 & -\sigma&&\\
                           &1       & -\sigma &  &\\
                           &        & \ddots       &\ddots     &\\
                           &        &         &   1 & -\sigma\end{matrix}\right]:
{\fake}_{n+1} X=\bigoplus_{p=1}^{n+1} X_p\to {\fake}_n X
$$
Note $R$ is not a functor on $\pro$-$Sup\ay$ but only on $(Sup\ay)^\N$. 
We write $r:X\to RX$ for the composite of the canonical maps $X\to``\prod$''$X\to RX$; $r$
is a natural transformation of functors of inverse systems of complexes, and a cofibration and weak 
equivalence of pro-complexes.
Now let $f:A\to B$
be as above. We may assume $f$ has a level representative $\{f_n:A_n\to B_n\}_n$. Since
the latter is a map of inverse systems, we may apply $R$ to it. 
Consider the pro-supercomplex $N=B\oplus RB[1]\oplus RA$ with boundary map 
$$
\left[\begin{matrix}\partial & 0 &0 \\
      -r& -\partial & Rf\\
       0 & 0 & \partial\end{matrix}\right]
$$ 
Put $j:=(f,0,r):A\to N$, and let $q:N\to B$ be the projection. It is clear 
that $f=qj$, that $q$ and $j$ are a respectively a fibration and a cofibration, and that
both are chain maps. Furthermore we have a commutative diagram 
$$
\xymatrix{C[1]\ar@{ >->}[r]&N\ar@{->>}[r]&RA\\
           &&A\ar[ul]^j\ar[u]^r}
$$
where the row is a cofibration sequence, and $C$ is the mapping cone of $B\to RB$. As $C$ is
weakly contractible, $N\fib RA$ is a weak equivalence. Since also $A\to RA$ is one, it follows
that $j$ is a weak equivalence.
\end{proof}
\begin{propo}\label{M1prop}
Let 
\begin{equation}\label{M1}
\xymatrix{A\ar@{ >->}[d]_i\ar[r]^\alpha& X\ar@{>>}[d]^\pi\\
          B\ar[r]_\beta\ar@{.>}[ur] &Y}
\end{equation}
be a commutative solid arrow diagram in $\pro$-$Sup\ay$, where $i$ is a cofibration and weak equivalence and $\pi$ is a
 fibration. Then the dotted arrow exists and makes it commute.
\end{propo}
\begin{proof}
Because $i$ is a cofibration there is a $p'\in\Hom_{\Z/2}(B,``\prod$''$A)$
such that $p' i:A\to ``\prod$''$A$ is the canonical cofibration. Because $\pi$ is a fibration,
there is a homogeneous (but not necessarily chain) map $j:Y\to X$ such that $\pi\circ j=1_Y$. Without
loss of generality we may assume all of the following

\begin{itemize}
\item $\alpha$, $\beta$, $i$ and $\pi$ have level representative maps such that for each $n$
the diagram
$$
\xymatrix{A_n\ar[r]^{\alpha_n}\ar[d]_{i_n}&X_n\ar[d]^{\pi_n}\\
          B_n\ar[r]_{\beta_n}& Y_n}
$$
commutes.

\item $j$ and $p'$ have representatives $p'_n:B_{n+1}\to A_n$, $j_n:Y_{n+1}\to X_n$, with
$p'_n  i_{n+1}=\sigma$, $\pi_n  j_n=\sigma$, $j_{n-1}\sigma=\sigma j_n$.

\item Put $Q_n=\coker i_n$, $K_n=\ker\pi_n$. There are sequences of homogeneous maps 
$ h_n\in\hom_{\Z/2}(Q_{n+1},Q_n[1])$, $s_n\in\hom_{\Z/2}(K_n,K_{n+1})$ such that 
$[\partial,h_n]=\sigma$, $\sigma^2s_n=\sigma$.
\end{itemize}
Write $p:B\to Q$ for the projection. Then $i'=\sigma-ip'$ satisfies
$i'i=0$ and thus induces a homogeneous (but possibly not chain) map $Q\to ``\prod$''$B$, which 
we still call $i'$; it is  represented by a sequence of maps $i':Q_{n+2}\to B_n$.  Put 
$$
I':=[\partial,i'h]
$$
Then $I'$ is represented by a sequence of chain maps $I':Q_{n+2}\to B_n$, and we have 
\begin{equation}\label{pI'}
p I'=\sigma^2
\end{equation}
Put $P'=\sigma^2-I'p$; as $pP'=0$, we regard $P'$ as a sequence of maps $B_{n+2}\to A_n$. 
By construction
\begin{equation}\label{P'}
I'p+iP'=\sigma^2,\qquad P'i=\sigma^2
\end{equation}
Put
$$
\phi':=\sigma\alpha P'\sigma+[\partial,j\beta I'h]p:B_{n+4}\to X_n
$$
One checks that
\begin{equation}\label{phi'}
[\partial,\phi']=0,\qquad \phi'i=\alpha\sigma^4,\qquad \pi\phi'=\beta\sigma^4
\end{equation}
In other words $\phi'$ represents a chain map $B\to ``\prod$''$X$ making the following diagram
commute
\begin{equation}\label{diafake}
\xymatrix{A\ar[d]_i\ar[r]^\beta &X \ar[r]&\fake X\ar[d]^{\fake\pi}\\
          B\ar[r]_\beta\ar[urr]^{\phi'}&Y\ar[r]&\fake Y}
\end{equation}
Put
$$
\psi_1=\phi'_1:B_{\phi'(1)}\to X_1
$$
Next assume $n\ge 2$ and recursively that $\psi_{n-1}:B_{\psi(n-1)}\to X_{n-1}$
has been defined, is a chain map and satisfies that the classes of $\psi_{n-1}i_{\phi'(n-1)}$
and $\pi_{n-1}\psi_{n-1}$ in $\colim_m\hom(A_m,X_{n-1})$ and $\colim_m\hom(B_m,Y_{n-1})$
agree respectively with those of $\alpha_n$ and $\beta_n$. Then by \eqref{diafake} there exists
a representative 
$$
\phi'_n:B_{\phi'(n)}\to X_n
$$
of the $n$-th coordinate of 
$\phi'$ in 
$$\prod_m\colim_r\hom(B_r,X_m)=\hom(B,\fake X)$$
such that $\phi'(n)>\psi(n-1)$ and such that for
$$
g_{n-1}:=\psi_{n-1}\sigma^{\phi'(n)-\psi(n-1)}-\sigma\phi'_n
$$
we have 
$$
\pi_{n-1} g_{n-1}=0
$$
Hence we may regard $g_{n-1}$ as a map $B_{\phi'(n)}\to K_{n-1}$. Put
$$
\psi_n:=\phi'_n\sigma^3+[\partial,sg_{n-1}I'hp]
$$
Then $\psi_n$ is a chain map and $\psi_ni_{\psi(n)}$ and $\pi_n\psi_n$ have respectively the
same class as $\alpha_n$ and $\beta_n$ in $\colim_m\hom(A_m,X_{n-1})$ and $\colim_m\hom(B_m,Y_{n-1})$. 
Furthermore, one checks, using the assumptions on $s$ and
$h$ together with \eqref{pI'} and \eqref{P'}, that 
$$
\sigma^2\psi_n=\sigma\psi_{n-1}\sigma^{\phi'(n)-\psi(n-1)+3}
$$
Put 
$$
\phi_n=\sigma\psi_n
$$
Then $\phi$ commutes with $\sigma$ and therefore represents a map $B\to X$; one checks furthermore
that this map makes \eqref{M1} commute.

\end{proof}

\begin{coro}\label{fibequiv}
Let $\pi:X\defor Y$ be a fibration and weak equivalence in $\pro$-$Sup\ay$. Then there exists a 
chain map $g:Y\to X$ such that $\pi g=1_Y$ and such that $g\pi$ is homotopic to $1_X$.
\end{coro}
\begin{proof}
Assume first that $Y=0$. Then $X$ is fibrant and weakly contractible, and we have to show it is 
contractible. Let $X\cofi CX$ be the inclusion
into the mapping cone of the identity. Consider the solid arrow diagram
$$
\xymatrix{X\ar[r]^{1}\ar@{ >->}[d]&X\ar@{->>}[d]\\
          CX\ar[r]\ar@{.>}[ur]& 0}  
$$
Because $X\to 0$ is a weak equivalence, so is $X\to CX$. Thus the dotted arrow exists by 
\ref{M1prop}. It follows that $X$ is contractible. Next we prove the general case. Let $\pi:X\defor Y$ be
as in the corollary, $K=\ker \pi$, $j:Y\to X$ a homogeneous map with $\pi j=1_Y$. By what we have just 
proved, there is a degree one map $h:K\to K$ with $[\partial, h]=1$. One checks that 
$g=j+h[\partial,j]$ has the required properties. 
\end{proof}

\begin{rema}
In \cite{gros} a closed model category structure for $\pro$-simplicial sets was introduced. 
In {\it loc. cit.}
cofibrations are defined as those pro-maps which can be represented, modulo isomorphism of maps
of pro-simplicial sets,
by an inverse system of cofibrations.  On the other hand it is not hard to see, using \ref{chri}, that
a map in $\pro-Sup\ay$ is a cofibration in our sense $\iff$ it is isomorphic, {\it as a map of 
pro-$\Z/2$-graded objects,} to a pro-map represented by an inverse system of cofibrations of 
$\Z/2$-graded objects, that is, of split monomorphisms.
\end{rema}

\section{Algebras}\label{alg}

\subsection{Algebras}\label{salg}
Let $(\ay,\ij)$ be as in section \ref{cplx}. Give $\ay$ the Waldhausen category structure where 
cofibrations are as defined in section \ref{cplx} and where the weak equivalences are the isomorphisms.
A {\it tensor product} in $\ay$ is a bifunctor
$$
\otimes:\ay\times\ay\to \ay
$$
such that the following two conditions are fulfilled. 
\begin{itemize}
\item For every $A\in\ay$,
$A\otimes$ and $\otimes A$ are additive exact functors of Waldhausen categories which preserve 
all cokernels.

\item $(\ay,\otimes)$ is a symmetric, strict monoidal category in the sense of \cite{mac}.
\end{itemize}
If $\otimes$ is a tensor product, we write
$$
t_{A,B}:A\otimes B\wiso B\otimes A
$$
for the symmetry isomorphism, and $k$ for the unit of $\otimes$. 
By an {\it algebra} with respect to $\otimes$
we understand a not necessarily unital monoid in the monoidal category $(\ay,\otimes)$. Explicitly
an algebra is an object $A\in\ay$ together with a map
$$
\mu:A\otimes A\to A
$$ 
(called {\it multiplication}) satisfying the usual associativity condition; no map $k\to A$ is assumed to exist. Whenever
such a map is given and satisfies the axiom for a unit, we say that 
$A$ is unital. If $A$ is an algebra, we write $\tilde{A}$ for $A\oplus k$ equipped with the 
multiplication
$$
\xymatrix{(\mu+1_A+1_A\oplus1_k):\tilde{A}\otimes\tilde{A}=A\otimes A\oplus A\oplus A\oplus k
\ar[r]&              \tilde{A}}
$$
Note $\tilde{A}$ is unital, the unit being given by the inclusion $k\to\tilde{A}$. The notions
of left, right and two-sided modules and submodules are the obvious ones. 
In particular (two sided) ideals are defined. An ideal $I$ of $A$ is called a {\it cof-ideal}
if $I\cofi A$ is a cofibration.  

\subsection{The cyclic mixed complex of an algebra}\label{4.2}
Fix $(\ay,\ij,\otimes)$. Write $\alg\ay$ for the category of algebras in $\ay$. 
For $A\in\alg\ay$
define a mixed complex $\Omega A$ as follows. Just as for usual algebras, put
$$
\Omega^0A=A,\quad\Omega^nA=\tilde{A}\otimes A^{\otimes n}\quad (n\ge 1)
$$
Write $p:\tilde{A}\to A$, and $i:k\to \tilde{A}$ for the obvious projection
and inclusion maps. We omit writing the inclusion map $A\to \tilde{A}$. Put
\begin{align}
t_n:=&t_{A^{\otimes n},A}:{A}^{\otimes n+1}\to {A}^{\otimes n+1}\nonumber\\
\mu_i:=&1_{\tilde{A}}^{\otimes i-1}
\otimes\mu\otimes 1_{\tilde{A}}^{\otimes n- i-1}:\tilde{A}^{\otimes n+1}\to
\tilde{A}^{\otimes n}\quad (0\le i\le n-1)\label{moo}\\
\mu_n:=&\mu_0  t_n  (p\otimes 1^{\otimes n}):\Omega^n A\to\Omega^{n-1}A\nonumber
\end{align}
and set
$$
b:\Omega^nA\to\Omega^{n-1}A,\quad b:=\sum_{i=0}^n(-1)^i\mu_i
$$
Define a map 
$$B:\Omega^nA=A^{\otimes n}\oplus A^{\otimes n+1}\to \Omega^{n+1}A$$
as zero on  $A^{\otimes n}$ and as 
$$
B=(i\otimes 1^{\otimes n+1}) \sum_{i=0}^nt_n^i \text{ \quad on } A^{\otimes n+1}.
$$
One has to check that $\Omega A$ is a mixed complex; this follows from the associativity of $\mu$ and 
the fact that $(\ay,\otimes)$ verifies the axioms (\cite{mac}) of a 
symmetric strict monoidal category with twisting map $t$. We observe that the
direct summand 
$$
C_nA=A^{\otimes n+1}\subset \tilde{A}\otimes A^{\otimes n}=\Omega^nA
$$
is mapped by $b$ to $C_{n-1}A$. Thus $C(A)$ is a complex, and the inclusion
$i:C(A)\cofi \Omega^nA$ is a split monomorphism, whence a cofibration. The bar
complex  $C^{bar}A$ is defined as the cokernel of $i$ shifted by $-1$; thus
$$
C(A)\cofi \Omega A\to C^{bar}A[-1]
$$
is a cofibration sequence. 
Next let $K\cofi A$ be a cof-ideal. We
want to define relative chain complexes $C(K:A)$, $C^{bar}(K:A)$ and a relative
mixed complex $\Omega(K:A)$ in such a way that each of the following
\begin{align}
\Omega(K:A)\to \Omega A&\to \Omega\frac{A}{K}\label{cosi}\\
 C(K:A)\to C(A)&\to C(\frac{A}{K})\label{cosii}\\
 C^{bar}(K:A)\to C^{bar}(A)&\to C^{bar}(\frac{A}{K})\label{cosiii}
\end{align}
be a cofibration sequence. Put $\Omega^0(I:A)=C_0(I:A)=C_0^{bar}(I:A)=I$. To define
the higher terms of these complexes, proceed as follows. For each fixed $n\ge 1$, 
consider the set 
$[1]^{[n]}$ of all maps 
$$[n]=\{0,\dots,n\}\to [1]=\{0,1\}$$
Define a partial order by $f\le g\iff f(i)\le g(i)$ for all $i\in [n]$. 
Consider the 
functor $$\Ca_n:[1]^{[n]}\to \ay$$
which assigns to each map  $f:[n]\to [1]$ the object
$$\Ca_n(f):=Y_{f(0)}\otimes\dots\otimes Y_{f(n)}$$
where $Y_0:=K$, $Y_1:=A$, and sends $f\le g$ to the corresponding map in $\ay$
induced by $K\cofi A$. Consider the colimit of the restriction of $\Ca_n$ to the 
complement of the constant map $x\mapsto 1$ ($x\in[n]$)

$$
C_n(K:A):=\colim_{[1]^{[n]}\backslash\{1\}}\Ca_n
$$ 

There is an induced map  

\begin{equation}\label{in}
i_n:C_n(K:A)
\cofi \colim \Ca_n=C_nA
\end{equation}

The direct sum of $i_n$ and $i_{n-1}$ gives

\begin{equation}\label{i}
\CD
\Omega^n(K:A):=C_n(K:A)\oplus C_{n-1}(K:A)@>i_n\oplus i_{n-1}>>\Omega^nA
\endCD
\end{equation}

The formula for the boundary $b$ of $\Omega^nA$ defines a map 
$\Ca_n(f)\to \Omega^{n-1}(K:A)$ for every 
$$f\in [1]^{[n]}\backslash\{1\}
\cup[1]^{[n-1]}\backslash\{1\}$$ which
lies in $C_{n-1}A$ if $f\in [1]^{[n]}\backslash\{1\}$. 
One checks that these are compatible and thus define a
$$b:\Omega^n(K:A)\to\Omega^{n-1}(K:A)$$
which restricts to $C_n(K:A)\to C_{n-1}(K:A)$. 
The Connes operator $B:\Omega^n(K:A)\to\Omega^{n+1}(K:A)$ is defined similarly.
One checks that these choices of boundary maps make $C(K:A)$ and $\Omega(K:A)$ 
respectively a chain and a mixed 
complex and that $C(K:A)\to C(A)$ and $\Omega(K:A)\to\Omega A$ are compatible
with boundaries. By construction, $C(K:A)\to \Omega(K:A)$ is a degreewise split
monomorphism, whence a cofibration. The relative bar complex is defined by the 
cofibration sequence
\begin{equation}\label{cobar}
C(K:A)\cofi\Omega(K:A)\to C^{bar}(K:A)[-1]
\end{equation}
It remains to show that, as was announced, \eqref{cosi},\eqref{cosii} and
\eqref{cosiii} are cofibration
sequences. First we need two lemmas.

\begin{lem}\label{prerel1}
Let 
$$
\xymatrix{A\ar@{ >->}[r]^{i} & X\\
          C\ar@{ >->}[u]^{\alpha}\ar@{ >->}[r]_{\beta}&B\ar@{ >->}[u]_j}
$$
be a commutative diagram where all four maps are cofibrations. Assume the
induced map 
$$
\bar{j}:B/C\cofi X/A
$$
is a cofibration too. Then also the induced map
$$
A+B:=A\coprod_CB\cofi X
$$
is a cofibration.
\end{lem}
\begin{proof} Let $I$ be a basic object and $f\vee g:A\coprod_{C}B\to I$ a map. 
We must show that $f\vee g$ extends to $X$. Because $i$ is a cofibration, $f$ 
extends to 
a map $p:X\to I$. One checks that for $h=g-p  j$, one has $h \beta =0$. Write
$\bar{h}:B/C\to I$ for the induced map. Because $\bar{j}$ is a cofibration,
$\bar{h}$ extends to a map $\bar{q}:X/A\to I$. We write $q$ for the composite
of $\bar{q}$ with the projection $X\to X/A$. One checks that $p+q$ extends 
$f\vee g$.
\end{proof}
\begin{lem}\label{prerel2}
Let $S\cofi M$ and $T\cofi N$ be cofibrations. Consider the coproduct
$$S\otimes N+M\otimes T:=(S\otimes N)\coprod_{S\otimes T}(M\otimes T)$$
Then
\begin{equation}\label{stseq}
S\otimes N+M\otimes T\to M\otimes N\to\frac{M}{S}\otimes\frac{N}{T}
\end{equation}
is a cofibration sequence.
\end{lem}
\begin{proof}
Consider the commutative diagram
$$
\xymatrix{S\otimes N\ar[r] & M\otimes N\\
          S\otimes T\ar[u]\ar[r] &M\otimes T\ar[u]}
$$
It follows from the exactness of $\otimes$ that the induced map between the 
cokernels of the vertical arrows is
$$
(S\cofi M)\otimes \frac{N}{T}
$$
which is a cofibration, again by exactness of $\otimes$. By \ref{prerel1} the first map
of \eqref{stseq} is a cofibration; it remains to
show the second map is the cokernel
of the first. This is straightforward from the exactness of $\otimes$.
\end{proof}

\begin{prop}\label{cosindeed}
Let $K\cofi A$ be a cof-ideal. Then \eqref{cosi}, \eqref{cosii} and \eqref{cosiii} are 
cofibration sequences.
\end{prop} 
\begin{proof}
It suffices to show that, for $n\ge 1$ and  $i_n$ as in \eqref{in},
\begin{equation}\label{seqn}
\xymatrix{C_n(K:A)\ar@{ >->}[r]&C_nA\ar[r]&
C_n(\frac{A}{K})}
\end{equation}
is a cofibration sequence. For $n=1$ this is a particular case of 
\ref{prerel2}. Let $n\ge 1$ and assume by induction that \eqref{seqn}
is a cofibration sequence. One checks that 
$$
C_{n+1}(K:A)\cong C_n(K:A)\otimes A\coprod_{C_n(K:A)\otimes K}
C_nA\otimes K
$$
The inductive step now follows using \ref{prerel2}.
\end{proof}

\section{Wodzicki's excision theorem}\label{mariusz}
\begin{defi}
Let $K\in\alg\ay$, and $D$ one of the functorial complexes $C$, $C^{bar}$, $\Omega$, $\B\Omega$ of the 
previous sections. We say that $K$ is {\it $D$-excisive} if for every cofibration $K\cofi A\in\alg\ay$
making $K$ into a cof-ideal, the natural map $D(K)\to D(K:A)$ is a weak equivalence.
\end{defi}
\begin{them}\label{wodteo}
The following are equivalent for $K\in\alg\ay$.
\item{i)}$V\otimes C^{bar}K$ is weakly 
contractible for every $V\in\ay$.
\item{ii)} $K$ is $C$-excisive.
\item{iii)}$K$ is $C^{bar}$-excisive.
\item{iv)} $K$ is $\Omega$-excisive.
\item{v)} $K$ is $\B\Omega$-excisive.
\end{them}

\begin{proof}
For usual algebras this is Wodzicki's theorem \cite{wod}, thm. 3.1. We follow
the Gucciones' proof (\cite{gg}, thm. 4) and make the appropriate observations
as to what changes are needed. Thus what follows should be read with a copy
of \cite{gg} at hand. First of all note that purity is not an issue
in our case; by definition of cof-ideal and exactness of $\otimes$, the cofibration 
sequence
$$
K\cofi A\cofi A/K
$$
remains one upon tensoring with
any $V\in\ay$. The proof of \cite{gg} thm. 4 has two preliminaries, namely
theorem 2 and corollary 3. Theorem 2 of {\it
loc. cit.} holds with the same proof, replacing `quasi-isomorphism' by
`weak equivalence' everywhere, and noting that the filtration in the proof
is a filtration by cof-subobjects. Similar considerations apply to the
proof of corollary 3. In addition one must also note two other things. The
first is that, for 
$$C(A,B(p))[-p-1]:=(B(p)\otimes A^{\otimes *-p-1},b)$$
there is a cofibration sequence
$$
C(A,B(p))[-p-1]\cofi \tilde{F}^p\to\tilde{F}^{p+1}
$$
This, together with the argument in {\it loc. cit.} shows that 
$(\tilde{F}^p\to \tilde{F}^{p+1})\otimes V$ is a weak equivalence for all 
$p$ and all $V\in\ay$. The second thing to note is that, by \ref{weq}, 
this suffices to prove $(C(A,B)\to C(B))\otimes V$ is a weak equivalence.
Now to the proof of \cite{gg} thm. 4 itself. In the proof that 
i)$\Rightarrow$ii), replace $\ker\pi$ by $C(K:A)$, and exact sequence
by cofibration sequence. The proof of i)$\Rightarrow$ii) stays the same.
That ii) and iii) imply iv) is immediate from the cofibration sequence
\eqref{cobar}. In the proof of ii)$\Rightarrow$i) the algebra structure given to 
$A=V\oplus K$ can be defined in terms of arrows rather than elements in an
obvious way. Also note that, as $\otimes$ preserves $\oplus$, $C(K:A)$
equals the kernel of map of complexes induced by the projection $A\to V$, so
that the proof does apply to our case. That iii)$\Rightarrow$i) is analogous.
That  iv)$\Rightarrow$v) follows from \ref{Bexact}. The proof of \ref{Bexact} together
with the $SBI$ sequence give v)$\Rightarrow$iv).
To prove iv)$\Rightarrow$i) proceed as follows. For $A$, $V$ as in 
the proof of thm. 4 of \cite{gg},
the relative forms decompose into a sum of complexes
$$\Omega(K:A)=\Omega K\oplus\bigoplus_{p=1}^\infty\Omega_p(K:A)$$
Here $\Omega^n_p(K:A)=0$ for $p\ge n+1$, and for $p\le n+1$, it
 is the sum of all summands of tensor products of $n+1-p$ 
$K$'s and $p$ $V$'s in $A^{\otimes n+1}$ plus the sum of all tensor products
of $n-p$ $K$'s and $p$ $V$'s in $k\otimes A^{\otimes n}$, excepting $k\otimes V^{\otimes p}$ when $n=p$. Because $\Omega K\cofi \Omega(K:A)$ is a weak equivalence
by hypothesis, $\Omega_p(K:A)$ is weakly contractible for all $p\ge 1$. 
Note that $V\otimes C^{bar}K[-1]$ is a subcomplex of $\Omega_1(K:A)$ and
that the inclusion is a split monomorphism degreewise. Write $Q$ for the quotient
complex. Thus $\Omega_1(K:A)$ is isomorphic to the
mapping cone of a chain map $f:Q[+1]\to V\otimes C^{bar}K[-1]$. Because
$\Omega_1(K:A)$ is weakly contractible, $f$ is a weak equivalence. I claim
that $Q$ is contractible. The claim implies that $V\otimes C^{bar}K$ is 
weakly contractible, as we had to prove. Before proving the claim, we introduce
some notation. Put
$$
P_n=K\otimes V\otimes K^{n-1}\oplus K^{\otimes 2}\otimes V\otimes K^{\otimes n-2}
\oplus\dots\oplus K^{\otimes n}\otimes V
$$
Then
$$
Q_n=P_n\oplus P_{n-1}
$$
and the boundary map has the following matricial form
$$
b=\left[\begin{matrix} \partial &\alpha\\ 
                 0  & \partial'\end{matrix}\right]
$$
where $\alpha:P_{n-1}\to P_{n-1}$ is of the form
$$
\alpha=\left[\begin{matrix} 1 & \theta_1 &&&\\
                     0 &1 &\theta_2&&\\
                       &  &\ddots &\ddots&\\
                     &&&1&\theta_{n-2}\\
                     &&&&1\end{matrix}\right]
$$
In particular $\alpha$ is an isomorphism. It follows that $Q$ is contractible,
whence also weakly contractible. \end{proof}

\begin{defi} 
We say that $K$ is $H$-unital if it satisfies the equivalent
conditions of \ref{wodteo}.
\end{defi}

\begin{stass}
Starting here and for the remaining of this paper we assume every object of $\ay$
is relatively injective.
\end{stass}

\begin{rema}
The standing assumption is equivalent to saying that all cofibrations of $\ay$ are split injections.
According to the definitions of \ref{pro}, this means that the cofibrations in $\pro$-$\ay$ are those maps
$A'\to A$ which have the extension property \eqref{rlp} for all $V\in\ay$. Note not all of these are 
split.
\end{rema}

\section{Excision for ideals of quasi-free pro-algebras}\label{exkinfi}

We say --after \cite{cq}-- that $A\in\pro$-$\alg\ay$ is {\it quasi-free} if the canonical inclusion
$$
d\cup d:A\otimes A\hookrightarrow \tilde{A}\otimes A\otimes A=\Omega^2A
$$
considered as a Hochschild $2$-cochain, is a coboundary. A {\it fundamental $1$-cochain} is
a map $\phi:A\to \Omega^2A$ whose coboundary is $d\cup d$.

\begin{exam}
Let $A\in\alg\ay$. We want to define an analogue of the pro-algebra $TA/JA^\infty$ of \cite{cq} in our setting.
Note that as we do not assume infinite sums exist in $\ay$, we have no analogue for the tensor
algebra $TA$. Instead we consider the pro-object
$$
\Tau A:=\fake \Omega^{2*}A
$$
and equip it with a pro-algebra structure as follows. First we consider the map
$$
\xymatrix{\Omega^{2r}A\otimes\Omega^{2s}A\ar[r]^{\mu_{r,s}}&\Omega^{2(r+s)}A\oplus\Omega^{2(r+s+1)}A\\
\tilde{A}\otimes A^{2(r+s)}\oplus\tilde{A}\otimes A^{2(r+s)+1}
\ar@{=}[u]&\tilde{A}\otimes A^{2(r+s)}\oplus\tilde{A}\otimes A^{2(r+s+1)}\ar@{=}[u]}
$$
defined by the matrix
$$
\left[\begin{matrix} 1& b'\otimes 1^{\otimes 2s}\\
              0 &- p\otimes 1^{\otimes{2(r+s)+1}}\end{matrix}\right]
$$
Here $p:\tilde{A}\to A$ is as in \ref{4.2}. 
Because $\otimes$ commutes with finite sums, we have a well-defined map
$$
\mu_n:=\sum_{0\le r,s\le n}\mu_{r,s}:\Tau_nA\otimes\Tau_nA\to \Tau_nA
$$
compatible with the structure maps $\sigma$. The proof of the associativity of the Fedosov
product of \cite{cq1} for usual algebras, shows that $\mu_n$ makes $\Tau_nA$ into an algebra,
and $\Tau A$ into a pro-algebra. To define a fundamental $1$-cochain for $\Tau A$, proceed
as follows. Assume first $\ay=\C$-vectorspaces. Take the fundamental cochain
given in \cite{cq1} for the tensor algebra $TA$, and rewrite it in terms of noncommutative
forms and Fedosov product. The result is the following formula, where we have written
$\delta$ for the de Rham differential of $TA$, $\omega_i=dx_idy_i$, for $x_i,y_i\in A$
and $\circ$ for the Fedosov product.
\begin{align*}
\phi(a\omega_1\dots \omega_n)=&\delta a\delta(\omega_1\dots\omega_n)
+\sum_{i=1}^{n-1}a\omega_1\dots\omega_{i-1}\delta\omega_i\delta(\omega_{i+1}\dots\omega_n)\\
-&\sum_{i=1}^na\omega_1\dots\omega_{i-1}\delta x_i\delta(y_i\omega_{i+1}\dots\omega_n)\\
+&\sum_{i=1}^na\omega_1\dots\omega_{i-1}\delta(x_i\circ y_i)\delta(\omega_{i+1}\dots\omega_n)\\
-&\sum_{i=1}^n(a\omega_1\dots\omega_{i-1}\circ x_i)\delta y_i\delta(\omega_{i+1}\dots\omega_n)
\end{align*}
This gives an induced map
$$
\phi_n:\Tau_nA\to\Omega^2\Tau_{n-1}A
$$
commuting with $\sigma$. Note that, because the formula above is expressed in terms of the Fedosov product and the product of the
algebra $A$, it can be written entirely in terms of arrows, and therefore makes sense for algebras
on an arbitrary additive category $\ay$. Moreover as the proof of the identity
$$
\phi(\xi\circ\eta)=\phi(\xi)\eta+\xi\phi(\eta)+\delta \xi\delta\eta
$$
for algebras on $\ay=\C-Vect$ is a formal consequence of the associativity and distributivity of the 
product of $A$, it holds for arbitrary $\ay$, in its proper element-free formulation.
In particular, the pro-algebra $\Tau A$ is quasi-free for all $A\in\alg\ay$.
\end{exam}

\begin{nota}
The $n$-th power $R^n$ of an algebra $R\in\alg\ay$ is defined as
the image (i.e. the cokernel of the kernel) of the $n$-fold multiplication map $R^{\otimes n}\to R$.
One checks (using the fact that $\otimes$ preserves cokernels) that if $R$ is an ideal of an 
algebra $S$ then the same is true of $R^n$.

If $A\in\pro-\alg\ay$ then a {\it pro-ideal} of $A$ is a sequence of ideals $K_n\triangleleft A_n$
($n\ge 1$) which is compatible with the maps $A_n\to A_{n-1}$. We write 
$K^\infty=\{K^n_n:n\in\N\}$. The pro-ideal $K\triangleleft A$ is
a {\it cof-pro-ideal} if the inclusion $K\cofi A$ is a cofibration. 
\end{nota}

\begin{them}\label{xinfi}
Let $A\in\pro$-$\alg\ay$ be quasi-free and
$K\cofi A$ a cof-pro-ideal. Then $K^\infty$ is $H$-unital
and the inclusion $K^\infty\cofi K$ is a cofibration.
\end{them}

\begin{proof}
We shall show first that if $A$ is quasi-free and $L\cofi A$ is a cof-pro-ideal such that $L^2\to L$
is an isomorphism then $L$ is $H$-unital. Then we shall prove that, for $K$ as in the proposition,
the map $K^\infty\to K$ is a cofibration. This suffices to prove the proposition, since 
$(K^\infty)^2\cong K^\infty$. For $K$ as above, the map $\tilde{\mu}$ in the commutative diagram
below gives ``$\prod$''$K$ a left $A$-pro-module structure
\begin{equation}\label{alf0}
\xymatrix{A_n\otimes{\fake}_n K\ar[rr]^(.55){\tilde{\mu}}\ar[d]_{1}&&{\fake}_n K\\
\bigoplus_{p\le n}A_n\otimes K_p\ar[rr]^{\sum\sigma^{n-p}\otimes 1}&&\oplus_{p\le n}A_p\otimes 
K_p\ar[u]_{\sum\mu_p}}
\end{equation}
A right $A$-module structure is defined analogously; the two actions make 
``$\prod$''$K$ into an $\tilde{A}$-pro-bimodule and the canonical cofibration 
$\iota :K\cofi$``$\prod$''$K$ into a morphism
of such. Because $K_n$ is relatively injective for all $n$, ``$\prod$''$K$ is relatively
injective in $\pro$-$\ay$ (by \ref{equi-inj}). Thus there exists a map 
$e:A\to$``$\prod$''$K\in\pro$-$\ay$ such that $ej=\iota$ where
$j:K\cofi A$ is the inclusion. Let $\beta$ be the composite

\begin{align*}
\beta:K\overset{j}\to A\overset{\phi}\to
\tilde{A}\otimes A\otimes A\overset{1\otimes 1\otimes e}\longrightarrow&\\
&\ \ \tilde{A}\otimes A\otimes
\fake K
\overset{b'}\to\tilde{A}\otimes\fake K
\end{align*}
Here $\phi$ is a choice of fundamental $1$-cochain and $b'$ is the boundary of the bar complex
with coefficients in the $A$-pro-bimodule ``$\prod$''$K$. For $i$ as in \ref{4.2}, put
$$
\alpha=\beta+i\otimes\iota
$$
One checks that the diagram
\begin{equation}\label{alf1}
\xymatrix{A\otimes K\ar[r]^(.4){1\otimes\alpha}\ar[d]_\mu &A\otimes\tilde{A}\otimes
\fake K\ar[d]^{\mu\otimes 1}\\
K\ar[r]_(.4){\alpha} &\tilde{A}\otimes\fake K}
\end{equation}
commutes. It follows that there is a commutative square
$$
\xymatrix{K^2\ar[r]^(.3){\psi_0}\ar[d]&K\otimes\fake K\ar[d]\\
          K\ar[r]_(.3){\alpha}&\tilde{A}\otimes\fake K}
$$
One checks further that also
\begin{equation}\label{alf2}
\xymatrix{K\ar[r]^(.3){\alpha}\ar[dr]_\iota&\tilde{A}\otimes\fake K\ar[d]^{b'}\\
&\fake K}
\end{equation}
commutes. Set 
\begin{align*}
\psi:C^{bar}_*(K^2)\to C^{bar}_{*+1}(K,\fake K)&\\
\psi_n:=(-1)^n {j'}^{\otimes n-1}\otimes\psi_0:(K^2)^{\otimes n}\to &
 K^{\otimes n}\otimes\fake K
\end{align*}
Here ${j'}:K^2\hookrightarrow K$ is the inclusion. 
One checks, using \eqref{alf1} and \eqref{alf2}, that 
$$
b'\psi_n+\psi_{n-1}b'={j'}^{\otimes n-1}\otimes \iota j'\qquad (n\ge 0)
$$
Next we observe that if $V,W\in\pro$-$\ay$ then 
\begin{align*}
\theta:V\otimes\fake W\to&\fake (V\otimes W),\\
\theta_n=\sum_p\sigma^{n-p}\otimes 1:
\bigoplus_{p\le n}V_n\otimes W_p&\to\bigoplus_{p\le n}V_p\otimes W_p
\end{align*}
is a natural map of inverse systems and satisfies
$$
\theta (1_V\otimes \iota_W)=\iota_{V\otimes W}
$$
Moreover for $\tilde{\mu}$ as in \eqref{alf0}, we have a commutative diagram of inverse systems
$$
\xymatrix{K\otimes\fake K\ar[d]_{\tilde{\mu}}\ar[r]^\theta&\fake (K\otimes K)
\ar[dl]^{\fake \mu_K}\\
\fake K&\\}
$$
It follows that 
$$
\theta_*:=\theta(K^{\otimes *-1},K):C^{bar}(K,\fake K)\to (\fake (C^{bar}K),
\fake b')
$$
is a homomorphism of pro-complexes. Note also that
$$
\theta_*  (1^{\otimes *-1}\otimes \iota)=\iota^{\otimes *}
$$
Thus for 
$$
\zeta:=\theta \psi:C^{bar}(K^2)\to \fake (C^{bar}K)
$$
we have 
$$
(\fake b')\zeta_n+\zeta_{n-1}b'=\iota  {j'}^{\otimes n}
$$
for all $n$. In other words, $\zeta$ is a homotopy $\iota  {j'}\sim 0$. Let $V\in\pro$-$\ay$. Then
$$
\zeta_V=1_V\otimes\zeta
$$
is a homotopy $1_V\otimes \iota  {j'}\sim 0$. Thus $\zeta_V$ induces 
$$
V\otimes C^{bar}_{*-1}(K^2)\overset{\bar{\zeta_V}}{\longrightarrow}
\frac{ V\otimes\fake C_*^{bar}(K)}
{(1\otimes\fake b')(V\otimes 
\fake C^{bar}_{*+1}K)}
$$
Consider the map
$$
\bar{\theta}:\frac{ V\otimes\fake C_*^{bar}(K)}
{(1\otimes\fake b')(V\otimes 
\fake C^{bar}_{*+1}K)}
\to \fake \frac{V\otimes C_*^{bar}(K)}{(1\otimes b')V\otimes C_{*+1}^{bar}(K)}
$$
induced by $\theta(V,C_*^{bar}K)$. We have a commutative 
diagram
$$
\xymatrix{\frac{V\otimes C^{bar}_*K^2}{(1\otimes b')(V\otimes C^{bar}_{*+1}K^2)}\ar[r]^{b'}\ar[d]_{\iota
  {j'}}
&V\otimes C^{bar}_{*-1}(K^2)\ar[d]^{\bar{\zeta_V}}\\
\fake (\frac{V\otimes C^{bar}_*K}{(1\otimes b') (V\otimes C^{bar}_{*+1}K)})&
\frac{ V\otimes\fake C_*^{bar}(K)}{(1\otimes\fake b')(V\otimes \fake C^{bar}_{*+1}K)}
\ar[l]_{\bar{\theta}} }
$$
If ${j'}$ happens to be an isomorphism, this says that the top horizontal row in the diagram above is a
cofibration,
whence $V\otimes C_*^{bar}K$ is weakly contractible by lemma \ref{weq}.
It remains to show that $K^\infty\to K$ is a cofibration. Because $j:K\cofi A$ is a cofibration, there
exists a map $e_0:A\to``\prod K$'' such that $e_0  j =\iota$. For $p$ as in \ref{4.2}, consider the 
composite
\begin{multline*}
e_1:K\overset\alpha\to\tilde{A}\otimes K\overset{p\otimes 1}\longrightarrow A\otimes K\\
\overset{e_0\otimes 1}\longrightarrow (\fake K)\otimes K\overset{\theta}\to\fake (K\otimes K)\overset{\fake \mu}\longrightarrow \fake K^2
\end{multline*}
We have
$$
e_1\circ (K^2\to K)=\iota_{K^2}
$$
Thus $K^2\cofi K$ and $K^2\cofi A$ are cofibrations. This process can be repeated inductively, to
obtain, for each $n$, a map 
$$e_n:K^{2^n}\to\fake K^{2^{n+1}}$$
such that
$$
\xymatrix{K^{2^{n+1}}\ar[r]\ar@{ >->}[d]_{\iota_{K^{2^{n+1}}}}&K^{2^n}\ar[dl]^{e_n}\\
\fake K^{2^{n+1}}&}
$$
commutes. Moreover, because $``\prod$''$K^{2^{n+1}}$ is relatively injective and 
$\iota_n:=\iota_{K^{2^n}}$ is a
cofibration, there exists a map $h_{n}$ such that 
$$
\xymatrix{&\fake K^{2^{n+1}}\\
\fake K^{2^{n}}\ar@{.>}[ur]^{h_n}&K^{2^n}\ar@{ >->}[l]_(.4){\iota_n}\ar[u]^{e_n}}
$$          
commutes. We need to show that there exists a map $u$ making the following diagram commute
$$
\xymatrix{K^\infty\ar[r]^j\ar[d]_{\iota_\infty}&A\ar@{.>}[dl]^u\\
\fake K^\infty&}
$$
The existence of such a map is equivalent to the existence, for each $m\ge 1$, of a map
$$
A_{u(m)}\overset{u}\longrightarrow K_m^m
$$
such that, for some $u'(m)\ge u(m)$,
\begin{equation}\label{qpq}
\xymatrix{K^{u'(m)}_{u'(m)}\ar[r]^\sigma\ar[dr]_\sigma&K^{u(m)}_{u(m)}\ar[r] &A_{u(m)}\ar[dl]^u\\
&K_m^m&}
\end{equation}
commutes. Put 
$$v=v_n:=h_{n-1}\circ h_{n-2}\circ\dots\circ e_0:A\to\fake K^{2^{n}}$$
Then 
$$
\xymatrix{K^{2^n}\ar[r]\ar[d]&A\ar[dl]^v\\
\fake K^{2^n}&}
$$
commutes. Thus we can choose a representative  
$$
\hat{v}\in\prod_p\hom(A_{v(p)},K^{2^n}_p)
$$
of 
$$v\in\hom(A,\fake K^{2^n})=\prod_p\colim_q(A_q,K^{2^n}_p)$$
such that 
$$
\forall p\quad \hat{v}_p\circ (K^{2^n}_{v(p)}\to A_{v(p)})=\sigma^{v(p)-p}
$$
For each given $m$ as in \eqref{qpq}, choose $n$ such that $2^n\ge m$. Then (omitting
the exponents of the transfer maps)
$$
\xymatrix{K^{v_n(2^n)}_{v_n(2^n)}
\ar[r]^\sigma\ar[dr]_\sigma &K^{2^n}_{v_n(2^n)}\ar[r] & A_{v_n(2^n)}\ar[d]^{v_n}\\
&K_m^m&K_{2^n}^{2^n}\ar[l]_\sigma}
$$
commutes. Thus the composite map $u=\sigma  v_n$  
satisfies \eqref{qpq}.
\end{proof}

\section{Goodwillie's theorem}

Let $M=\{M^m\}_m$ be a pro-mixed complex. A {\it cof-filtration} of $M$ is a filtration
$M=\F_0\supset\F_1\supset\dots$ 
by sub-mixed complexes
such that $\F_{n+1}\cofi\F_n$ is a cofibration for all $n$. We call $\F$ {\it bounded} if
for every level $m$ there exists an $n_k$ such that $X^m\F_{n_k}^m=0$.

\begin{lemm}\label{aboveG}
Let $M$ be a pro-mixed complex and  $\F$ a bounded cof-filtration on $M$. Assume
further that there is a nondecreasing, nonstationary function $f:\N\to \N$ such that for 
every $n$
the inclusion $\F_n\cofi \F_1$ has a representative
 $$
\{\F_n^{f(m)}\to\F_1^m\}_m
$$

If  $X^\infty(\F_n/\F_{n+1})$ is weakly contractible for all $n\ge 1$, then
$X^\infty M\fib X^\infty (M/\F_1)$ is a weak equivalence.
\end{lemm}
\begin{proof}
Because $X^\infty$ preserves cofibration sequences (by \ref{xiexact}), it suffices to show that 
$X^\infty \F_1$ is weakly contractible. An induction argument shows that
the cofibration $i:X^\infty \F_n\to X^\infty \F_1$ is a weak equivalence for all $n$. 
For $f$ as in the lemma we have a commutative diagram 
\begin{equation}\label{abog1}
\xymatrix{\F_n\ar@{ >->}[r]\ar[d]& \F_1\ar[d]\\
           \F_n^{f(m)}\ar[r] &\F_1^m}
\end{equation}
whence also a commutative solid arrow diagram
\begin{equation}\label{abog2}
\xymatrix{X^\infty\F_n\ar@{ >->}[r]\ar[d] & X^\infty\F_1\ar@{.>}[dl]\ar[d]\\
          X^m\F^{f(m)}_n\ar[r]&X^m\F_1}
\end{equation}
Because $\F_n\cofi \F_1$ is a cofibration and weak equivalence and $X^m\F^{f(m)}_n$ is a constant
supercomplex, the diagonal arrow exists and makes the upper left triangle commute. It follows
that the lower right triangle commutes in $Ho\pro$-$Sup\ay$. But because $\F$ is bounded, for each 
fixed $m$ one
can choose $n$ so that  $\F_n^f(m)=0$. It follows that for each $m$ the map 
$X^\infty \F_1 \to X^\infty \F^m_1$ is nullhomotopic. Hence $X^\infty \F_1$ is weakly contractible.
\end{proof}
\begin{lemm}\label{aboveG2}
 Let $M,N\in\ay$, 
\begin{align*}
M=\F_0\supset\F_1&\supset\F_2\supset\dots\\
N=\G_0\supset\G_1&\supset\G_2\supset\dots
\end{align*}
cof-filtrations. For $1\le m\le n$ define inductively
$$
\sum_{i=0}^m\F_i\otimes\G_{n-i}:=(\sum_{i=0}^{m-1}\F_i\otimes\G_{n-i})\bigoplus_{\F_m\otimes \G_{n-m+1}}
\F_m\otimes \G_{n-m}
$$
and put $(\F\otimes\G)_n=\sum_{i=0}^n\F_i\otimes \G_{n-i}$. Then $\F\otimes \G$ is a cof-filtration
of $M\otimes N$ and 
$$\frac{(\F\otimes\G)_n}{(\F\otimes\G)_{n+1}}=\bigoplus_{i=0}^n\frac{\F_i}{\F_{i+1}}\otimes\frac{\G_{n-i}}
{\G_{n-i+1}}$$
\end{lemm}
\begin{proof}
To prove $\F\otimes\G$ is a cof-filtration we must show $(\F\otimes\G)_n\to M\otimes N$ is a cofibration
for all $n$. We shall prove by induction on $1\le m\le n$ that 
\begin{equation}\label{sumene}
\sum_{i=0}^m\F_i\otimes \G_{n-i}\to M\otimes N
\end{equation}
is a cofibration. For $m=0$ this map is 
$\F_0\otimes\G_n=M\otimes\G_n\cofi M\otimes N$ which is a cofibration because $\G_n\cofi N$ is
and $\otimes$ is exact. Assume that $1\le m\le n$ and that 
$\sum_{i=0}^{m-1}\F_i\otimes \G_{n-i} \to M\otimes N$ is a cofibration. To prove \eqref{sumene}
is a cofibration, it suffices --by lemma \ref{prerel1}-- to show that the top row of the commutative
diagram below is a cofibration 
\begin{equation}\label{coke}
\xymatrix{\F_m\otimes \frac{\G_{n-m}}{\G_{n-m+1}}\ar@{ >->}[r]\ar@{ >->}[d]&
 \frac{M\otimes N}{\sum_{i=0}^{m-1}\F_i\otimes \G_{n-i}}\ar@{->>}[d]\\
\F_m\otimes \frac{N}{\G_{n-m+1}}\ar@{ >->}[r] &M\otimes \frac{N}{\G_{n-m+1}}}
\end{equation}
But this follows from the fact that both the left vertical and lower horizontal arrows are.
The remaining assertion of the lemma is straightforward from the exactness of $\otimes$.
\end{proof}
\begin{them}\label{goo}
Let $\ay$ be $\Q$-linear, $A\in\pro$-$\alg\ay$, 
$$
A_m=\F_0A_m\supset \F_1A_m\supset\F_2A_m\supset\dots \quad (m\ge 1)
$$
a multiplicative filtration of $A_m$ compatible with the maps $A_{m+1}\to A_m$. Assume 
$\F_*A=\{\F_*A_m\}_m$ is cof-filtration and is such that, for each $m$, there exists an
$n_m$ such that $\F_{n_m} A_m=0$.
Then
$$
X^\infty(A)\defor X^\infty(A/\F_1A)
$$
is a weak equivalence.  
\end{them}
\begin{proof}
By \ref{aboveG2} and induction we get that the filtration $\F:=\F A$ induces on $\Omega A$ is a cof-filtration
such that $\F_r\Omega A/\F_{r+1}\Omega A$ is isomorphic to the pro-complex $\Omega_rG$ which 
at each level $m$ is the homogenous part of degree $r$ of the forms on the pro-graded algebra 
$G_m=\bigoplus_{r\ge 0}\F_rA_m/\F_{r+1}A_m$. 
By \ref{aboveG} it suffices to show that $X^\infty(\Omega_r G)$ is contractible
for $r>0$. The proof of \cite{lod}, 4.1.11 applied to the Euler
derivation shows that for $r>0$, each of the maps 
$X^n\Omega_rG_m\to X^{n-1}\Omega_rG_m$ 
is nullhomotopic. Thus $X^\infty\Omega_r G$ is weakly contractible, as we had to prove.
\end{proof}
\section{Cuntz-Quillen theorem}\label{seccq}
\begin{them}\label{cq}
Assume $\ay$ is $\Q$-linear.
If $A\in\pro-\alg\ay$ and $I\triangleleft A$ is a cof-pro-ideal, then $X^\infty I\to X^\infty(I:A)$ is
a weak equivalence.
\end{them}
\begin{proof}
Let $\Tau A:=\{\Tau_nA_n\}_n\fib \Tau A/I$ be the induced map, $K$ its kernel. At each level $n$, the 
map
\begin{equation}\label{tseqn}
K_n\cofi \Tau_nA_n
\end{equation}
is the sum for $p\le n$ of the maps
\begin{equation}\label{oseqn}
\Omega^{2p}(I_n:A_n)\to \Omega^{2p}A_n
\end{equation}
Because $I\cofi A$ is a cofibration, 
$$
\Omega^{2p}(I:A)\cofi \Omega^{2p}A
$$
is a cofibration by \ref{cosindeed}. Therefore for each $n$ there exist an $N$ and a map
$\Omega^{2p}A_N\to\Omega^{2p}(I_n:A_n)$ 
making the following diagram commute
\begin{equation}\label{soluno}
\xymatrix{\Omega^{2p}(I_N:A_N)\ar[r]\ar[d]&\Omega^{2p}A_N\ar[dl]\\
         \Omega^{2p}(I_n:A_n)&}
\end{equation}
Because for each level $n$, $K_n\to\Tau_nA$ involves only finitely many terms \eqref{oseqn},
\eqref{soluno} implies that $K\to \Tau A$ is a cofibration. 
Consider the descending filtrations of $K$, $\Tau A$ and $\Tau A/I$ induced by the decomposition
of \eqref{oseqn}. Each of these satisfies the hypothesis of \ref{goo}. It follows
that both vertical arrows in the commutative diagram
$$
\xymatrix{X^\infty K\ar[r]\ar[d]_{\wr}&X^\infty(K:\Tau A)\ar[d]^{\wr}\\
          X^\infty I\ar[r] &X^\infty(I:A)}
$$
are weak equivalences. Hence it suffices to show the top horizontal arrow is a weak equivalence.
Consider the commutative diagram
$$
\xymatrix{X^\infty K^\infty\ar[r]^{\sim}\ar[d]_{\wr}&X^\infty(K^\infty:\Tau A)\ar[dd]\\
 X^\infty (K^\infty:K)\ar[d]&\\         
X^\infty K\ar[r]& X^\infty(K:\Tau A)}
$$
The decorated maps are weak equivalences as follows from applying successively \ref{xinfi}, \ref{wodteo}, 
\ref{xiexact} and \ref{remxiexact}.
By \ref{xinfi}, the filtration $\F=\{\F^n\}_n$ of $K_n/K_n^\infty$ given by 
$\F^n_0=\F_1^n=K_n/K_n^n$, $\F_r^n= K_n^r/K_n^n$ $(r\le n)$ satisfies the hypothesis of 
\ref{goo}. Thus  $X^\infty(K^\infty:K)\to X^\infty K$ is a weak equivalence. Similarly,
$X^\infty(\Tau A/K^\infty)\to X^\infty\Tau(A/I)$ 
is a weak equivalence. From the commutative
diagram

$$
\xymatrix{
X^\infty(K^\infty:\Tau A)\ar[r]\ar[d]&X^\infty\Tau A\ar[r]\ar@{=}[d]&X^\infty(\Tau A/K^\infty)\ar[d]^\wr\\
X^\infty(K:\Tau A)\ar[r]&X^\infty\Tau A\ar[r]&X^\infty(\Tau A/I)}
$$
we get that also $X^\infty(K^\infty:\Tau A)\to X^\infty(K:\Tau A)$ is a weak equivalence. 
\end{proof}

\begin{defi}
Let $A,B\in\pro$-$\alg\ay$. The {\it bivariant periodic cyclic hypercohomology} of $(A,B)$
is  
$$
\HP^n(A,B)=\hom_{Ho\pro-Sup\ay}(X^\infty A,X^\infty B[n])
$$ 
Here $\hom$ is taken in the homotopy category of $\pro$-$Sup\ay$. 
\end{defi}
\begin{rema}
It follows from \cite{qui}, I.1, Thm.1, that for arbitrary pro- algebras $A,B$, 
$$
\HP^n(A,B)=H^n\Hom_{\Z/2}(X^\infty A,\widehat{X}^\infty B)
$$
where the $\widehat{\ \ }$ indicates an arbitrary choice of fibrant resolution 
$$
X^\infty B\overset\sim\cofi\widehat{X}^\infty B
$$
in $\pro$-$Sup\ay$.
For example one can take $\widehat{X}^\infty B=RX^\infty B$, where $R$ is as in the proof of \ref{clmod}. 
In the particular case when $B$ is an algebra,
$X^\infty B$ is fibrant already, whence 
$$
\HP^n(A,B)=H^n\Hom_{\Z/2}(X^\infty A,X^\infty B)=:HP^n(A,B)
$$
is the usual periodic cyclic theory.
\end{rema}
\begin{coro}\label{6tes}
Let $A,B\in\pro$-$\alg\ay$, $I\triangleleft A$ a cof-ideal. Then there are exact sequences
$$
\xymatrix{\HP^0(A/I,B)\ar[r]&\HP^0(A,B)\ar[r]&\HP^0(I,B)\ar[d]\\
          \HP^1(I,B)\ar[u]&\HP^1(A,B)\ar[l]&\HP^1(A/I,B)\ar[l]\\
          \HP^0(B,A/I)\ar[d]&\HP^0(B,A)\ar[l]&\ar[l]\HP^0(B,I)\\
          \HP^0(B,I)\ar[r]&\HP^0(B,A)\ar[r]&\HP^1(B,A/I)\ar[u]}
$$
\end{coro}
\begin{proof}
The theorem shows that 
$$X^\infty I\to X^\infty A\to X^\infty (A/I)$$
is both a fibration and a cofibration sequence (\cite{qui}, I.3) in $Ho\pro-$ $Sup\ay$.
The corollary then follows from \cite{qui}, I.3, Prop. 4 and 4'.
\end{proof}

\begin{ack}
Important discussions for this research were carried out during visits of each of the authors to
the other's institution, as well as to The Abdus Salam ICTP. Funding for these visits was provided 
by ICTP, CONCYTEC, IMCA and the 
University of Buenos Aires. We are duly grateful to all of them. The second author wishes to
thank the PUCP for reducing his teaching load so that he could carry out this investigation.
\end{ack}

\end{document}